\documentclass[preprint,12pt]{article}

\usepackage{amsthm,amsfonts,amssymb,mathrsfs,graphics,natbib}

\usepackage{latexsym}
\usepackage{amssymb}
\usepackage{amsmath}
\usepackage[english]{babel}
\usepackage{fontenc}
\usepackage{times}
\usepackage[all]{xy}
\usepackage{epsfig}
\usepackage{amsthm}
\usepackage{inputenc}
\usepackage{setspace}
\usepackage{yfonts}
\usepackage{lscape}
\usepackage{natbib}
\usepackage{url}
\usepackage{fixltx2e}
\usepackage{alltt}
\usepackage{graphicx}
\usepackage{mathrsfs}
\usepackage{wrapfig}

\graphicspath{{immagini//}}

\usepackage[all]{xy}
\usepackage{setspace}
\usepackage{lscape}

\newcommand{\C}{\mathbb{C}}
\newcommand{\CP}{\mathbb{CP}}

\newcommand{\R}{\mathbb{R}}
\newcommand{\Gr}{\mathbb{G}\mathrm{r}}
\newcommand{\E}{\Lambda}

\renewcommand{\theenumii}{\arabic{enumii}}
\def\s#1.{\stepcounter{enumii}\theenumi.\theenumii\ \textsl{#1.}\\}

\newcommand{\y}{\\[0pt]}

\newcommand\mb{\smallbreak}
\newcommand\n{\noindent}
\newcommand\ip{\raise1pt\hbox{\large$\lrcorner$}\,}

\newcommand{\g}{\mathfrak{g}}

\newcommand{\ot}{\otimes}

\newcommand{\we}{\wedge}
\newcommand{\ba}{\begin{array}}\newcommand{\ea}{\end{array}}
\newcommand{\rf}[1]{(\ref{#1})}\def\lie#1({\mathfrak{#1}(}
\newcommand{\ft}[2]{\hbox{$\textstyle\frac{#1}{#2}$}}
\newcommand{\e}{\mathrm{e}}
\newcommand{\nit}[1]{\mb\n\textit{#1.}}

\def\a{\ol a}

\def\,{\kern2pt}
\def\.{\,\cdot\,}

\theoremstyle{plain}
\newtheorem{teor}{Theorem}
\newtheorem{prop}[teor]{Proposition}
\newtheorem{lemma}[teor]{Lemma}
\newtheorem{corol}[teor]{Corollary}
\theoremstyle{definition}
\newtheorem{defi}[teor]{Definition}
\theoremstyle{remark}
\newtheorem{rem}[teor]{Remark}

\begin{document}

\title{Intrinsic torsion classes of \\ Riemannian structures}
\author{Georgi Mihaylov}
\maketitle
\begin{abstract}\noindent This article introduces
  the notion of intrinsic torsion varieties associated to
  $G$-structures on a parallelizable Riemannian manifold. As an
  illustration, the intrinsic torsion varieties of orthogonal almost
  product structures are analysed on the Iwasawa manifold.
\end{abstract}

\section*{Introduction}

Determining integrability of a geometrical structure is a
fundamental problem in differential geometry. A suitable formulation
relates the equivalence problem for $G$-structures with the theory
of $k$-jets. In this context, the intrinsic torsion of a
$G$-structure is a first-order obstruction to integrability. In
certain circumstances, vanishing of the intrinsic torsion is not
only a necessary but also a sufficient integrability condition. The
celebrated Darboux and Newlander--Nirenberg theorems assert exactly
this for almost symplectic and almost complex structures
respectively. On the other hand, in the Riemannian situation,
non-vanishing intrinsic torsion impedes a reduction of the holonomy
group.\smallbreak

Symmetry properties of the intrinsic torsion tensor $\tau$ have been
exploited in various cases for classifying manifolds with
$G$-structure. The prototype classification is given by Gray and
Hervella in [7], and regards almost Hermitian and almost
symplectic manifolds. Later in [11], Naveira applied
analogous approach to the classification of almost product
manifolds, Cabrera and Swann studied the case of
quaternion-Hermitian manifolds [4] etc. Applications of a
Gray--Hervella type classification intervene in the construction of
specific $G_2$-structures on 7-manifolds (see [5]).
Constructions of symplectic and other structures of special classes
are described in [6] and in other articles by Salamon,
Fino, Chiossi and others.\smallbreak

In this article we consider a manifold $M$ with a fixed Riemanian
metric and subgroup $G$ of $O(N)$. Different $G$-reductions of the
Riemannian structure are parametrized by a suitable space.\smallbreak

A favorable context to start our analysis is the case in which $M$
is a nilmanifold. Nilmanifolds are parallelizable in a natural way,
their geometry reflects properties of the underlying nilpotent
algebra in terms of invariant tensors. An invariant $G$-structure
becomes a point of the homogeneous space $O(N)/G$, and this space
contains the subsets parameterizing special $G$-structures as real
subvarieties. These are the \emph{intrinsic torsion varieties}
 (ITV's) associated to the problem.\smallbreak

In the first section, we introduce the precise definition of ITV's.
In the rest of the article we analyse a specific case, namely the
classification of orthogonal almost product structures on the
Iwasawa manifold. Similar analysis on almost Hermitian structures
(using different methods) has been developed in [1]. The
method described in the present article exploits essentially basic
definitions and so appears more general. The rich geometry of the
ITV's highlighted in these cases (by rather simple calculations)
motivates the research of a general theory which explains the
variation of the intrinsic torsion of distinguished reductions of a
$G$-structure.

\section{Intrinsic torsion varieties}

The restriction of the canonical soldering form $\theta$ on the
frame bundle $LM$ of a smooth manifold to a $G$-structure $P$ allows
one to view the tangent space in each point of $M$ as a
representation $\rho_P$ of $G$ on $\R^n$. The covariant derivative
$\nabla_\phi \theta$ associated to a connection $\phi$ on $LM$
defines a horizontal torsion form $\Theta$, with values in
$\mathrm{Hom}(\R^n\wedge\R^n,\R^n)$. The \emph{structure function}
on $P$ represents, the component of $\Theta$ independent of the
choice of connection compatible with the structure, i.e. the image
of $\Theta$ in the space\begin{equation}
\mathscr{W}=\mathrm{Hom}(\R^n\wedge\R^n,\R^n)/\alpha(\mathrm{Hom}(\R^n,\g)).
\label{coker}
\end{equation}
Here $\alpha$ is the skew-symmetrization in $\mathrm{Hom}(\R^n,\g)$,
and $\g$ is the Lie algebra of $G$. The structure function defines
over $M$ a section $\tau$ of a vector bundle with fibre
\eqref{coker} called the \emph{intrinsic torsion of $P$}. The
variation of $\Theta$ along a fibre of $P$ defines a representation
$\rho_P$ on $\mathscr{W}$. If at a point $p\in P$ $\Theta$ belongs
to a $\rho_P$-submodule $\mathscr U$, in the projected point $\tau$
belongs to the same subspace. Equivalently, we say that $\tau$ maps
to zero in the quotient $\mathscr{W}/\mathscr{U}$. Adopting the
expression coined in [9]:

\begin{defi}
A \emph{null-torsion structure} on a manifold is a $G$-structure for
which $\tau$ belongs to a proper $G$-submodule $\mathscr{U}$ of
$\mathscr{W}$.
\end{defi}

The space \eqref{coker} decomposes into a direct sum of
$G$-irreducible subspaces $\mathscr{W}_i$. A classification of
$G$-structures based on criteria whereby $\tau$ belongs to a
specific subset of these components appears ``natural'' as $\Theta$
is basic tensorial form and ``geometrically relevant'' as the
obstructions determined by the intrinsic torsion reduce to a
specific subspace of $\mathscr{W}$.\smallbreak

Consider a fixed Riemannian manifold $M$. For $O(N)$, the space
\eqref{coker} is trivial, and the consequent vanishing of the
structure function explains the existence of a canonical
(Levi-Civita) torsion-free connection $\phi_{LC}$. Let $G$ be a
closed subgroup of $O(N)$. A reduction to a $G$-structure $P_\xi$
can be defined by detecting $G$ as the stabilizer $Satb(\xi)$ of a
tensor $\xi$. The parameter space of such reductions can be
point-wise identified with the $O(N)$-orbit $\mathscr{O}_\xi$ of
$\xi$, so a $G$-structure on $M$ is a smooth section of a bundle
with fibre $\mathscr{O}\cong O(N)/G$.\smallbreak

The torsion of the reduced structure comes from the restriction of
$\theta$ to $P_\xi$ eliminating the dependence of $\nabla\theta$ on
connections compatible with $P_\xi$. In general, the Levi-Civita
connection does not reduce to $P_\xi$, so $\tau_\xi$ does not
vanish. This incompatibility arises from the non-zero projection
onto $\mathfrak{g}^\bot_\xi$ (the orthogonal complement of $\g_\xi$
inside the vertical space of $P$) of a vector horizontal with
respect to $\phi_{LC}$. In other words, $\tau_\xi$ \emph{depends on
the mutual positions of the horizontal spaces} of $\phi_{LC}$ and
the connections compatible with $P_\xi$. For this reason in the
Riemannian case:
\begin{equation}
\mathscr{W}\cong\R^n\otimes\g^\bot_\xi. \label{iosperp}
\end{equation}

\noindent \textbf{Fact:} Reductions of an $O(N)$-structure may
belong to different null-torsion classes.\smallbreak \noindent This
fact has been remarked in [1, 9]. We can consider a
``local problem'' i.e. given a small open neighborhood of a point on
$M$ (using a single coordinate chart domain), determine the
mechanism which causes the class-variation of $\tau_\xi$ between
distinguished reductions of $P$.\smallbreak

Consider a Riemannian reduction $P_\xi$. A well-known fact is that:
\begin{equation}\mathrm{Stab}(g\xi)=g\,\,\mathrm{Stab}(\xi)\,\,g^{-1},\,\,\,\,\, g\in O(N).  \label{varstab}\end{equation}
\noindent Given a decomposition $\mathscr{W}=\oplus_i\mathscr{W}_i$,
the action of $g\in O(N)$ transforms it in $\oplus_i g
\mathscr{W}_i$. In fact given $h\in \mathrm{Stab}(\xi)$, suppose
$\tau_\xi\in \mathscr{W}_i$. As
$(ghg^{-1})(g\mathscr{W}_i)=g(h\mathscr{W}_i)$, $ghg^{-1}$ acts on
$g\mathscr{W}_i$ as $h$ acts on $\mathscr{W}$. So $P_\xi$ and
$P_{g\xi}$ belong to the same class if:
\begin{equation}\label{gttg} \tau_\xi\in \oplus_i\mathscr{W}_i\,\,\,\,\, and \,\,\,\,\, \tau_{g\xi}\in \oplus_i g \mathscr{W}_i \end{equation} \noindent The mere $O(N)$-action on the fibres on $P$ keeps $\tau_\xi$ in the same submodule.\smallbreak

The Levi-Civita connection form restricted to $P_\xi$ can be written
as $$\phi_{LC}=\phi_\xi+\psi_\xi.$$ \noindent Here $\phi_\xi$ is a
connection compatible with $P_\xi$ and $\psi_\xi$ represents the
projection on the specific orthogonal complement $\g^\bot_\xi$.
Along a fibre of $P$ the form $\phi_{LC}$ varies via the adjoint
representation of $O(N)$ and the variation of the splitting
$\g_\xi\oplus \g^\bot_\xi$ follows the rule \rf{varstab}. This
implies that $Ad_g\phi_\xi$ defines point-wise a connection
compatible with $P_{g\xi}$. As we just project on the subspaces of
the splitting, the variation of $\psi_\xi$ along the vertical
directions transversal to the fibres of $P_\xi$ is determined to
first order by the $O(N)$-action. So apparently the variation of all
elements (the representation on $\mathscr{W}$ and the obstruction
$\psi$) is determined by the $O(N)$-action. But this cannot justify
the null-torsion class variation. This means that a point-wise
first-order $G$-structure theory approach cannot explain this
phenomenon. The mechanism, which induces the ``jumps'' of the
reductions from one class to another should be expressed in terms a
higher-order prolongation of $P$. This consideration appears quite
natural in the context of the theory of prolongations of
$G$-structures (see for example [13]). In fact the
splitting of the tangent space in each point of $P$ determines a
$G$-structure over the total space of the $P$ i.e. a prolongation.
In [9] has been already underlined the relation between
the intrinsic torsion and the second fundamental form of the
embedding $P_\xi\subset P$, which contains more information.\smallbreak

A problem closely related to the local-one appears when $M$ is a
parallelizable manifold. If we declare a global section $s$ of the
frame bundle to be orthonormal, we define a Riemannian metric on
$M$. The bundle of orthonormal frames is a product $M\times O(N)$.
We call a $G$-structure \emph{invariant} if (relative to the
parallelization $s$) $\xi$ is constant. Invariant $G$-structures are
parametrized by a single \emph{classifying orbit} $\mathscr{O}$.

\begin{defi} Intrinsic torsion variety (ITV) is a subset of $\mathscr {O}=O(N)/G$ composed by invariant reductions which belong to a given null-torsion class.
\end{defi}

\noindent\textbf{Problem:} Analyse the geometry of the Intrinsic
Torsion Varieties of invariant structures on a parallelizable
Riemannian manifold.\smallbreak

Each point $\xi$ in $\mathscr O$ defines an invariant reduction
$P_\xi$, and determines both, the intrinsic torsion $\tau_\xi$ and
the $G$-representation on the space \eqref{coker}. The intrinsic
torsion can be regarded as a map from $\mathscr{O}$ to
$\mathscr{W}$. The $G$-representation on $\mathscr{W}$ varies over
$\mathscr{O}$ by \rf{varstab}. Equation \rf{gttg} suggests a special
class of ITV's, namely constant torsion varieties (CTV's) can be
defined by the condition:
\begin{equation}\tau_{g\xi}=g\tau_\xi\label{ctveq}.\end{equation}

The $O(N)$-action on the tensor components can be expressed in terms
of polynomials, so $\xi \longrightarrow g\xi $ and
$\tau_\xi\longrightarrow g\tau_\xi$ are polynomial functions. This
means that CTV's can be described as zero-sets of polynomial
functions on $\mathscr{O}$ and so CTV's are \emph{varieties} in some
general algebraic sense. This observation, involving invariant
polynomial functions, suggests that tools of Geometric Invariant
Theory can be applied to the study of ITV's.  The key point of the
general description of ITV's is the variation of the intrinsic
torsion ``transversal'' to the $O(N)$ action.
\begin{rem}\label{ctvrem}We expect that an ITV is fibration of algebraic CTV's over a ``cross section'' in the set of $O(N)$ orbits in $\mathscr{W}$.\end{rem}

\section{Almost product structures}

Consider  a smooth section $\mathfrak{F}$ of the tensor bundle
$T^*M\otimes TM$ over a differentiable manifold. The ``value'' in
$p$ of $\mathfrak{F}$ is an endomorphism of $T_pM$. If
$\mathfrak{F}$ is stabilized by a subgroup $G$ of $GL(N,\R)$, it
defines a $G$-structure on $M$. The action of $\mathfrak{F}_p$ can
be analysed in terms of its eigenspaces which give rise to vector
distributions in $TM$. The integrability of these distributions
characterizes $M$. Almost complex, $f$-structures [2, 14],
mixed structures [9, 10] etc. are relevant examples.\smallbreak

The almost product structures (APS's) introduced by Naveira in [11] can be
included in this setup. An APS corresponds to a
splitting of the tangent space at each point of $M$ into a direct
sum of vertical and horizontal subspaces
$T_pM=\mathcal{V}\oplus\mathcal{H}$. Denoting by
$\mathfrak{v},\mathfrak{h}$ the projections on $\mathcal{V}$ and
$\mathcal{H}$, an APS can be defined stabilizing:
\begin{equation}\label{1way}\mathfrak{P}=\mathfrak{v}-\mathfrak{h}\,\,\,\in\,\,\,
T^*M\otimes TM.\end{equation}

An orthogonal almost product structure (OPS) is a reduction of a
Riemannian structure with group $O(p)\times O(q)$, $\mathcal{V}$ and
$\mathcal{H}$ are orthogonal subspaces. The same splitting regards
the cotangent space and the summands $\mathcal{V}$ and
$\mathcal{H}$ can be identified with their dual spaces via the
metric. The parameter space of these reductions is the Grassmannian
of $p$-planes in $\R^N$.\smallbreak

The intrinsic torsion of an OPS is given by the Levi-Civita
derivative $\nabla\mathfrak{P}$, equivalently by the Levi-Civita
derivative of the covariant 2-tensor $\Phi=g(\mathfrak{P}\,\,\,
,\,\,\,)$. As $$\mathfrak{so}(N)\ \cong\ \E^2(\mathcal{V}\oplus
\mathcal{H})=\E^2\mathcal{V}\oplus
\E^2\mathcal{H}\y\oplus(\mathcal{V}\we \mathcal{H})
\cong\mathfrak{so}(p)\oplus\mathfrak{so}(q)\oplus(\mathcal{V}\otimes
\mathcal{H}),$$ \noindent the space of intrinsic torsion
\rf{iosperp} is: \begin{equation} \R^N\otimes \mathfrak{g}^\perp
\cong(\mathcal{V}\oplus\mathcal{H})\otimes(\mathcal{V}\otimes
  \mathcal{H}).\label{compare}\end{equation}
OPS's has been classified by Naveira in [11] by means of
the following $O(p)\times O(q)$-irreducible components of the
intrinsic torsion space
\rf{compare}:\begin{equation}\begin{array}{lll}
    \mathscr{W}_1=\Lambda^2 \mathcal{V}\otimes \mathcal{H},\,\, & \mathscr{W}_2=S_0^2\mathcal{V}\otimes\mathcal{H},\,\, & \mathscr{W}_3=Tr(\mathcal{V}\otimes\mathcal{V})\otimes \mathcal{H},\\
    \mathscr{W}_4=\Lambda^2 \mathcal{H}\otimes \mathcal{V},\,\, & \mathscr{W}_5=S_0^2 \mathcal{H}\otimes \mathcal{V}, \,\, & \mathscr{W}_6=Tr(\mathcal{H}\otimes\mathcal{H})\otimes \mathcal{V}.
    \end{array}\label{classes}\end{equation}\noindent The null-torsion conditions detect thirty-six
classes of OPS's obtained by summing two or more of these subspaces
and exchanging ``vertical'' and ``horizontal''.\smallbreak

Let $v$ denote a section of the distribution $\mathcal{V}$. The
Levi-Civita connection on $P$ defines a connection on the vector
bundle $\mathcal{V}$ via the projections of $\nabla v$ on
$\mathcal{V}$. Vice versa the incompatibility of the Levi-Civita
connection with an OPS is determined by the projection of $\nabla v$
on $\mathcal{H}$. More precisely:

\begin{teor}\label{opercrit}The intrinsic torsion of an OPS is completely
determined by the
maps:$$\tau_v:\mathcal{V}\to\mathcal{H}\otimes\mathcal{H},
\,\,\,\,\,\,
\tau_h:\mathcal{H}\to\mathcal{V}\otimes\mathcal{V}$$\noindent
defined by the orthogonal projections of $\nabla v$ on
$\mathcal{H}\otimes\mathcal{H}$ and $\nabla h$ onto
$\mathcal{V}\otimes\mathcal{V}$.\end{teor} \nit{Proof} We take an
orthonormal basis of forms $\{v^i,h^j\}$, and express $\Phi$
as:$$\Phi=\sum_{i=1}^p v^i\otimes v^i -\sum_{j=1}^q h^j\otimes
h^j$$\noindent Thus $\nabla_w\Phi$ with $w\in TM$ is determined by
terms of the form:
\begin{equation}\label{nablap}\sum_{i=1}^{p} (\nabla_w v^i)\otimes
v^i-\sum_{j=1}^q (\nabla_w h^i)\otimes h^i\end {equation} The tensor
$\Phi$ restricted to $\mathcal{V}$ and $\mathcal H$ is proportional
to the metric. This implies $\nabla_x\Phi(y,z)=0$ for each $x\in TM$
if $y$ and $z$ are both vertical or horizontal. The non-zero
contributions of $\nabla_w v^i$ and $\nabla_w h^i$ to \rf{nablap}
belong respectively to $\mathcal{H}$ and $\mathcal{V}$ and
\rf{nablap} belongs to the space
$(\mathcal{V}\oplus\mathcal{H})\mathcal{V}\mathcal{H}$. So setting $
w=\sum_{k=1}^p a_kv_i+\sum_{l=1}^q b_l h_l$, we re-write \rf{nablap}
as:
$$ \begin{array}{l}  \sum_{i=1}^p\sum_{k=1}^p a_k(\nabla_{v_k}v^i)\otimes v^i + \sum_{i=1}^p\sum_{l=1}^q b_l(\nabla_{h_l}v^i)\otimes v^i-\\
 \,\,\,\,\,\,\,\,\,\,\,\,\,\,\,\,\,\,\,\,\,\,\,\,\,\,\,\,\,\,\,\,\,\,\,\,\,\,\sum_{j=1}^q\sum_{k=1}^pa_i(\nabla_{v_i}h^j)\otimes h^j-\sum_{j=1}^q\sum_{l=1}^q b_l(\nabla_{h_l}h^j)\otimes h^j.\end{array}$$ Compare now the first term, say $I$, (which belongs to
$\mathcal{V}\otimes\mathcal{H}\otimes\mathcal{V}$)
 to the third term $III$ (which belongs to
 $\mathcal{V}\otimes\mathcal{V}\otimes\mathcal{H}$). Let us write $I$
 (respectively $III$) in terms of the basis $\{h^j\}$ (respectively $\{v^i\}$):$$\begin{array}{rl}
 I  &= \sum_{i=1}^p\sum_{k=1}^p a_k\sum_{n=1}^q(\nabla_{v_k}v^i,h^n)h^n\otimes v^i\\
 -III&= \sum_{j=1}^q\sum_{k=1}^p  a_i\sum_{m=1}^p(\nabla_{v_i}h^j,v^m)v^m\otimes h^j
\end{array}$$
As $g(v,h)=0$ and $\nabla g=0$ we have
$(\nabla_wv,h)=-(v,\nabla_wh)$, so: $$\begin{array}{rl}
 I  &= \sum_{i=1}^p\sum_{k=1}^p a_k\sum_{n=1}^q(h^n,\nabla_{v_k}v^i)h^n\otimes v^i\\
 III&= \sum_{j=1}^q\sum_{k=1}^p  a_i\sum_{m=1}^p(h^j,\nabla_{v_i}v^m)v^m\otimes h^j
\end{array}$$
\noindent as the order of the summations and the names of the
indices are irrelevant, these terms are equivalent from the point of
view of the representation. The same consideration applied to the
second and the fourth term completes the proof.\qed\smallbreak

As the vertical space $\mathfrak{so}(n)\cong\Lambda^2 T_p$, a
suitable starting point for the general analysis of ITV's might be
the Riemannian $G$-reductions defined by stabilizing a $2$-form. In
the case of an (oriented) Riemannian manifold ``contraction'' of a
$2$-form with the metric determines skew-symmetric endomorphisms
$\mathfrak{F}_p$. We can define an OPS  by setting
$\mathcal{V}=\mathrm{Im}\mathfrak{F}_p,\mathcal{H}=\mathrm{ker}\mathfrak{F}_p$.
In particular a \emph{simple} $2$-form (one which with a suitable
basis can be written as $e\wedge f$ with $e$ and $f$ $1$-forms)
determines a reduction of the structure group to $SO(2)\times
O(N-2)$.\smallbreak

In [1] the Gray--Hervella classes of invariant Hermitian
structures on real six-dimensional nilmanifolds were embedded in
$\CP ^3$. The null-torsion conditions has been re-formulated in
terms of differential, Hodge star operator and wedge product. This
approach fails in our case. Consider $SO(2)\times SO(4)$ structures
on $\R^6$.\begin{equation}\label{dalpha}d\alpha\in\ba{rcl}
\E^3(\R^6)^*\ =\ \E^3(\mathcal{V}\oplus
 \mathcal{H}) &\cong&
\E^2\mathcal{V}\,\mathcal{H}\oplus
\mathcal{V}\,\E^2\mathcal{H}\oplus \E^3\mathcal{H}
\ea,\end{equation} \noindent where $\alpha$ is the simple $2$-form.
A ``complementary'' simple 4-form is $\beta=\ast \alpha$. Then
\begin{equation}\label{dbeta} d\beta\ \in\ \E^5(\R^6)\,
\cong\mathcal{V}\oplus \mathcal{H}.\end{equation}

\noindent Comparing \rf{compare} and \rf{classes} to \rf{dalpha} and
\rf{dbeta}, we conclude that, between them, $d\alpha$ and $d\beta$
determine only \emph{four} of the \emph{six} irreducible components
of intrinsic torsion. In [3] a subgroup $G$ of $GL(n, R)$
is called admissible if $G$ is the largest subgroup that fixes the
space $\Lambda^\ast(\R^N)^G$ of invariant forms on $\R^n$, $G$ is
strongly admissible if, on $G$-structures, the closedness of the
invariant forms is (point wise) equivalent to the vanishing of the
intrinsic torsion.

\section{Invariant SO(4)$\times$SO(2) structures on the Iwasawa manifold}

\noindent The Iwasawa manifold $\mathscr I$ is the set of right
cosets of the complex Heisenberg group over a lattice. So $\mathscr
I$ is a complex manifold. As a real manifold $\mathscr I$ admits a
left-invariant basis of real 1-forms $e^i$ such that:
\begin{equation}
de^1=0,\,\,\,de^2=0,\,\,\,de^3=0,\,\,\,de^4=0,\,\,\,
de^5=e^{13}+e^{42},\,\,\,de^6=e^{14}+e^{23}\label{oneforms}\end{equation}
\noindent where $e^{ij}=e^i\wedge e^j$. The Heisenberg group is
2-step nilpotent. The real geometry of $\mathscr{I}$ reflects the
natural decomposition of $\g^*=\mathscr{K}\oplus \langle
e^5,e^6\rangle$, where:$$\mathscr{K}=\ker(d|\g^*).$$\noindent The
map $\pi:G_H\to \C^2$ defined for $g\in G_H$ by $g\to (z^1,z^2)$
induces a fibration $\mathscr{I}\to T^4$ with fibre $T^2$. The space
tangent to the fibre is $\mathscr{K}^\bot=\langle e_5,e_6\rangle$.
The fact that $\mathscr{I}$ is a $2$-torus bundle over a $4$-torus
reflects a general property of nilmanifolds, their geometry can be
understood in terms of tower of toric fibrations.\smallbreak

An invariant tensor on $\mathscr{I}$ can be expressed in terms of
the bases $\e^i$ and $e_i$ with constant coefficients. Declaring
$e^i$ are orthonormal, we determine an invariant Riemannian metric
on $\mathscr{I}$. An invariant $SO(2)\times SO(4)$-structure is
determined by an invariant simple $2$-form. These OPS's are
parametrized by the Grassmannian $\Gr_2(\R^6)$. We single out a
generic $2$-plane in $\R^6$ by a couple of orthogonal unit $1$-forms
$v^1$ and $v^2$, or equivalently by the simple $2$-form
$\omega=v^1\wedge v^2$. Obviously such a parametrization of the
Grassmannian is  $SO(2)$-overabundant.\smallbreak

We can apply the Frobenius Theorem for establishing the
integrability of the distributions of an $SO(2)\times
SO(4)$-structure on $\mathscr{I}$:

\begin{prop} Given on $\mathscr{I}$ an orthonormal basis of vector fields
  $\{h_i\}$,

1) the distribution $\langle h_1,...,h_4\rangle$ is integrable if
and only if $d(h^{56})=0$

2) the distribution $\langle h_1, h_2\rangle$ is integrable if and
only if $d(h^{3456})=0$.\label{frob1}
\end{prop}

\nit{Proof} A tangent distribution $V\subset TM$ is integrable if
and only if its annihilator $V^\circ\subset\Omega^1M$ generates a
differential ideal. Consider $V^\circ=\langle v^1,v^2\rangle$.
\[dv^1=v^1\wedge\alpha+v^2\wedge \beta,\qquad dv^2=v^1\wedge\alpha'+v^2\wedge
\beta',\] and so \begin{equation}d(v^1\wedge v^2)=(v^1\wedge
v^2)\wedge\gamma\label{diffid}\end{equation} for some 1-form
$\gamma$. The nilpotency condition forces $d(v_1\wedge v_2)=0$. We
can set:
\begin{equation}\omega=v^1\wedge v^2=f^1\wedge f^2+f^3\wedge e^5+f^4\wedge
e^6+ae^{56}\label{f1234}\end{equation} \noindent with
$f^1,...,f^4\in \mathscr{K}$, then \rf{diffid}
becomes:\begin{equation} d\omega=f^3\wedge e^{13}+f^3\wedge
e^{42}+f^4\wedge e^{14}+f^4\wedge
e^{23}+a(e^{514}+e^{523}+e^{613}+e^{642})\label{diffid2}\end{equation}$$=f^{12}\wedge\gamma+f^3\wedge
e^5\wedge \gamma +f^4\wedge e^6\wedge
\gamma+ae^{56}\wedge\gamma$$\smallbreak The lack of $f^{12}$ term on
the left side forces $\gamma$ to be a linear combination of $f^1$
and $f^2$. If $a\ne0$, the same consideration applied to the
$e^{56}$-term shows that $\gamma$ is a linear combination of $e^5$
and $e^6$, which is a contradiction. If $a=0$,
$$f^3\wedge(e^{13}+e^{42})=f^3\wedge e^5\wedge\gamma\,\,\,\,\,,\,\,\,\,\,
f^4\wedge(e^{14}+e^{23})=f^3\wedge e^6\wedge\gamma$$ which again
contradicts the above condition. We conclude that the equality
\rf{diffid2} cannot be realized except in the case
$v_1,v_2\in\mathscr{K}$ when both sides are zero. On a nilmanifold,
analogous argument works for any subspace with annihilator whose
generators define a simple $k$-form with $k\geq 2$. This includes
the second item.\qed\smallbreak

Obviously from \rf{diffid2} follows that for a simple $2$-form
$\omega$:
\begin{equation}d\omega=0\,\,\,\,\Leftrightarrow\,\,\,\,\omega \in\E^2\mathscr{K}.\label{closedf12}\end{equation}

In this light we determine the set of OPS's characterized by closed
$4$-form complementary to the 2-form in (\ref{f1234}). The simple
$4$-form in question equals
\begin{equation}\label{simple4}
\ast\omega = g^{34}\wedge e^{56}-g^{234}\wedge e^6-g^{134}\wedge e^5
+g^{1234},
\end{equation}
in terms of suitable 1-forms $g^i$. The last term gives no
contribution to $d\ast\omega$, so there is no constraint on the
$e^{56}$-component of $\omega$. Furthermore we have
$$d(g^{234}\wedge
e^6)=g^{234}\wedge(e^{14}+e^{23})=0,\,\,\,\,\,\,\,d(g^{134}\wedge
e^5)=g^{134}\wedge(e^{13}+e^{42})=0,$$ \noindent these terms being
$5$-forms on a four-dimensional space. Finally,
$$d(g^{34}\wedge e^{56})=g^{34}\wedge(e^{13}+e^{42})\wedge e^6
- g^{34}\wedge e^5\wedge(e^{14}+e^{23}).$$ The vanishing of
$d\ast\omega$ is now seen to be equivalent to the
equations:\begin{equation}g^{34}\wedge\beta_2=0,\qquad g^{34}
  \wedge\beta_3=0,\label{beta1beta2}\end{equation}where $\beta_i$ is the following basis of the space $\Lambda^2_+\mathscr{K}$ of self-dual 2-forms on $\mathscr{K}$.
\begin{equation}\beta_1=e^{12}+e^{34},\,\,\,\,\,
\beta_2=e^{13}+e^{42},\,\,\,\,\,\beta_3=e^{14}+e^{23}\label{betas}\end{equation}

Any simple $2$-form $\gamma\in\Lambda^2(\R^4)^*$ can be written as
$\gamma=\sigma+\tau$ with $\sigma\in\Lambda_-^2$ and
$\tau\in\Lambda_+^2$ of equal norm. If $\gamma=g^{34}$, Equation
\rf{beta1beta2} forces $\sigma =\pm\beta_1$, so we obtain:
\begin{equation}\label{g34bt}
g^{34}=\pm\beta_1+\tau,
\end{equation} where $\tau$ is any unit element of $\Lambda^2_-$. Such 2-forms in
$\Lambda^2\mathscr{K}$ represent oriented real planes
$\left<v,J_1v\right>$ or $\left<v,-J_1v\right>$ invariant under the
action of the complex structure $J_1$ associated to the $2$-form
$\beta_1$. The entire set of such planes is parametrized by two
disjoint 2-spheres, determined by the sign in (\ref{g34bt}) and the
choice of $\tau$.\smallbreak

Returning to (\ref{f1234}) a simple form $\omega$ compatible with
the correct choice of $g^{34}$, has orthogonal projection on the
span $\langle\beta_2,\beta_3\rangle$ equal to zero. The set of such
forms is thus the real six-dimensional ``slice'' determined by the
kernel of this projection. We see that $\ast(f^{12})$ must be
$J_1$-invariant so $f^{12}$ is itself $J_1$ invariant.

\begin{rem}\label{condstar} In conclusion the 2-forms characterized by a closed complenetary 4-form determine a 2-plane with a $J_1$-invariant projection on $\Gr_2(\mathscr{K})$. We will frequently refer to this property as ``condition $(\ast)$''.\end{rem}

Combining Proposition~\ref{frob1}, \rf{closedf12} and \rf{g34bt} we
conclude:

\begin{prop} The set $\mathscr{V}$ of invariant $SO(2)\times SO(4)$-reductions of the standard Riemannian structure on $\mathscr{I}$ such that $[\mathcal{V},\mathcal{V}]\subseteq\mathcal{V}$ is the real 6-dimensional set of 2-planes satisfying condition ($\ast$). The set $\mathscr{H}$ characterized by $[\mathcal{H},\mathcal{H}]\subseteq\mathcal{H}$ is $\Gr_2(\mathscr{K})$. The set $\mathscr{D}=\mathscr{H}\cap\mathscr{V}$ is a disjoint union of two 2-spheres. \label{integrability}\end{prop}

The integrability of both $\mathscr{H}$ and $\mathscr{V}$ does
\emph{not} imply that the OPS is parallel. The holonomy of the
$\mathscr{I}$ cannot reduce and $SO(2)\times SO(4)$ is not in
Berger's list.

\section{ITV's of invariant OPS's on the Iwasawa manifold}

In this section we analyse the intrinsic torsion of invariant OPS's
on the Iwasawa manifold, exploiting the global orthonormal basis
$e^i$ \rf{oneforms}.  We denote by $\nabla$ the Levi-Civita
derivative, by $J_i$ the orthogonal almost complex structures
induced on $\mathscr{K}$ by the $2$-forms $\beta_i$ in \eqref{betas}
and by $e^i\odot e^j=\ft12(e^i\ot e^j+e^j\ot e^i)$ the symmetric
product. The following result is easily verified using the methods
of [12].

\begin{prop} Let $v$ denote a $1$-form on $\mathscr{I}$, $v_K$ its orthogonal projection on $\mathscr{K}$ and
$ae^5+be^6$ the one on $\mathscr{K}^\bot$. For the standard
Riemannian structure on $\mathscr{I}$:\begin{equation}2\nabla
v=(J_2v_K)\odot e^5+(J_3v_K)\odot
e^6+a\beta_2+b\beta_3.\end{equation}
\label{cov1forms}\end{prop}\noindent In view of Proposition
\ref{cov1forms} we characterize an OPS on $\mathscr{I}$ by the
intersections of $\mathcal{V}$ with $\mathscr{K}$ and
$\mathscr{K}^\bot$ i.e. by the $(\dim P\cap\mathscr{K}, \dim
P\cap\mathscr{K}^\bot)$. We denote by $\mathbf{T}_{(   \dim
P\cap\mathscr{K},\dim P\cap\mathscr{K}^\bot)}$ the subsets of
$\Gr_2(\R^6)$ of OPS's of a fixed type. Given an orthonormal basis
$v^1,v^2$ of $\mathcal{V}$, we put in evidence the components
\begin{equation} v^1=\mathbf f+\mathbf d,\,\,\, v^2=\mathbf
g+\mathbf
e,\,\,\,\,\,\,\mathbf{f},\mathbf{g}\in\mathscr{K},\,\,\mathbf{d},\mathbf{e}\in
\mathscr{K}^\bot\label{v1v2} \end{equation} \noindent Such a basis
can be completed to a basis of $T^*\mathscr{I}$ in two steps. First
we take an orthonormal basis $\{v^1,v^2,v^3,v^4\}$ of the span
$\langle \mathbf{f},\mathbf{d},\mathbf{g},\mathbf{e}\rangle$. Then
we get an orthonormal basis of $\mathcal{H}$, by choosing $v^5$,
$v^6$ in $\langle
\mathbf{f},\mathbf{g}\rangle^\bot\subset\mathscr{K}$. The precise
expressions and normalizing coefficients of the $v^i$-s are
irrelevant and therefore omitted.  Contrarily the form of the basis
with emphasis on the components in  $\mathscr{K}$ and
$\mathscr{K}^\bot$ is relevant.\medbreak

 \noindent \begin{tabular}{|c|c|c|}\hline
     Type  & Orthonormal basis $\mathcal{V\oplus\mathcal{H}}$ & $\mathbf{T}_{nm}$\\  \hline
     (0,2) & $\langle e^5,\pm e^6\rangle\oplus\langle v^1,v^2,v^3,v^4\rangle$ with $v_1,v_2,v_3,v_4\in\mathscr{K}$ & two points\\ \hline
     (2,0) & $\langle v^1,v^2\rangle\oplus\langle v^3,v^4,e^5,e^6\rangle$ with $v_1,v_2,v_3,v_4\in\mathscr{K}$ & $\Gr_2(\mathscr{K})$ \\ \hline
     (1,1) & $\langle \mathbf{f},\mathbf{e}\rangle\oplus\langle J_1\mathbf{f}, J_2 \mathbf{f}, J_3\mathbf{f},\mathbf{d}\rangle$ with $\mathbf{f}\in \mathscr{K}$, $\mathbf{e}\bot \mathbf{d}\in\mathscr{K}^\bot$ & $(S^3\times S^1)/\mathbb{Z}^2$\\ \hline
     (0,1) & $\langle \mathbf{f}+\mathbf{e},\mathbf{d}\rangle\oplus\langle J_1\mathbf{f}, J_2 \mathbf{f}, J_3\mathbf{f},\mathbf{f}-\mathbf{e}\rangle$ with $\mathbf{e}\bot \mathbf{d}$ & $S^5/\mathbb{Z}^2$\\ \hline
     (1,0) & $\langle \mathbf{f},\mathbf{g}+\mathbf{e}\rangle\oplus\langle \mathbf{g}-\mathbf{e}, v_3, v_4,d\rangle$ with $\mathbf{e}\bot \mathbf{d}$, $\mathbf{f}\bot \mathbf{g}$ & $S^4$ bundle over $S^3$ $/\mathbb{Z}^2$\\ \hline
     (0,0) & $\langle v^1,v^2\rangle\oplus\langle v^3,v^4,v^5,v^6\rangle$ as in \rf{v1v2}...& $\Gr_2(\R^6)$\\  \hline
  \end{tabular}

  \noindent Table 1: Types of $SO(2)\times SO(4)$-structures and compatible bases.\medbreak

\noindent In case $(0,1)$ we denote by $\mathbf{f}-\mathbf{e}$ a
unit form orthogonal to $\mathbf{f}+\mathbf{e}$ in the span $\langle
\mathbf{f},\mathbf{e}\rangle$, analogously $\mathbf{g}-\mathbf{e}$
in the case $(1,0)$. The topological description of each of these
sets is easily obtained by consecutive choices of the relevant
elements of the bases in the second column. The first two cases are
obvious. In the case $(1,1)$ we can choose $\mathbf{f}\in S^3$ and
$\mathbf{e}\in S^1$ and such a choice eliminates the $SO(2)$
overabundance of bases of the plane. The only residual ambiguity
arises from the fact that $(\mathbf{f},\mathbf{g})$ and
$(-\mathbf{f},-\mathbf{g})$ represent the same point in
$\Gr_2{\R^6}$. This is true for all the remaining types of OPS's
with fixed non-zero intersection with $\mathscr{K}$ or
$\mathscr{K}^\bot$. In the case $(1,0)$ the basis is fixed by the
decomposition $v_1=\mathbf{f}+\mathbf{e}$ of any unit $1$-form,
$\mathbf{d}$ is then determined automatically. In the case $(1,0)$
we choose a unit form $f\in S^3\subset \mathscr{K}$, then a unit
form in the span of $\langle J_1\mathbf{f}, J_2\mathbf{f},
J_3\mathbf{f}, e^5, e^6\rangle$.\smallbreak

\begin{rem}\label{type} Given $\pi\in \Gr_2(\mathscr{K})$ spanned by $v_1, v_2$, let $\varphi:\pi\longrightarrow \mathscr{K}^\bot$ be a linear map, then $v_1+\varphi(v_1)$ and  $v_2+\varphi(v_2)$
span a lifted plane $\pi^*\in \Gr_2(\R^6)$. We can identify
$\Gr_2(\R^6)$ with the total space of a fibre bundle $\mathcal{G}$
over $\Gr_2(\mathscr{K})$ with fibre
$\mathrm{Hom}(\pi,\mathscr{K})\cong \R^4$.  In these terms we
describe the subset $\mathscr{V}$. It can be identified with the
total space of a bundle $\mathcal{J}$ with fibre $\R^4$ over
$\mathscr{D}\cong\CP^1\sqcup \CP^1$. The ``compactification'' of
$\mathcal{G}$  and $\mathcal{J}$ is done performing in each fibre
the limits:\smallbreak

-- the zero section of $\mathcal{G}$ and $\mathcal{J}$ gives
$\mathbf{T}_{20}$ and $\CP^1\sqcup \CP^1\subset \mathbf{T}_{02}$;

-- for $\|\varphi(v)\|\longrightarrow \infty$ in each direction in
$\pi$ we reach $\mathbf{T}_{02}$;

-- for $\|\varphi(v)\|\longrightarrow\infty$ on a one-dimensional
subspace of $\pi$ we reach $\mathbf{T}_{01}$;

-- for $\|\varphi(v)\|\longrightarrow 0$ on a one-dimensional
subspace of $\pi$ we reach $\mathbf{T}_{10}$;

-- for $\|\varphi(v)\|\longrightarrow 0$,
$\|\varphi(v)\|\longrightarrow \infty$ on distinguished subspaces we
get $\mathbf{T}_{11}$.\smallbreak

\noindent OPS's in $\mathbf{T}_{02}$, $\mathbf{T}_{11}$ and
$\mathbf{T}_{01}$ satisfy trivially the condition $(\ast)$ (recall
Remark \ref{condstar}). \end{rem}\smallbreak

\noindent \textbf{Notation:} We call $\mathscr{C}_i$ the
null-torsion class complementary to $\mathscr{W}_i$ (recall
\rf{classes}) i. e. the $i$-th is the only missing component. Then
$\mathscr{W}_{ijk...}:=\mathscr{W}_i\oplus\mathscr{W}_j\oplus\mathscr{W}_k\oplus...$
and
$\mathscr{C}_{ijk...}:=\mathscr{C}_i\cap\mathscr{C}_j\cap\mathscr{C}_k\cap...\,\,\,$.
We call $\mathbf{W}_{ijk...},\,\,\, \mathbf{C}_{ijk...}\subset
\Gr_2(\R^6)$ the ITV's of points belonging respectively to the
null-torsion classes $\mathscr{W}_{ijk...}$ and
$\mathscr{C}_{ijk...}$. \smallbreak

\begin{teor} The ITV's of invariant $SO(2)\times SO(4)$ reductions of the standard Riemannian structure on the Iwasawa manifold are given by the following table:

\noindent \begin{tabular}{|l|l|l|}
  \hline ITV                                           &  topological description                & geometric properties               \\ \hline \hline
  $\mathbf{W}_4$                                       &  two points                             & totally geodesic horizontal foliations \\ \hline
  $\mathbf{W}_5\,\,\,\,\,\,\,\,\,(\equiv\mathscr{D})$  &  $\CP^1\sqcup \CP^1$                    & totally geodesic horizontal foliations \\ \hline
  $\mathbf{W}_{15}\,\,\,\,\,(\equiv\mathscr{H})$       &  $\Gr_2(\R^4)$                          & vertical foliations                    \\ \hline
  $\mathbf{W}_{345}$                                   &  $\CP^2\sqcup \CP^2$                    & horizontal foliations                  \\ \hline
  $\mathbf{W}_{245}$                                   &  $S^3\times S^1\times S^1$              & horizontal foliations                  \\ \hline
  $\mathbf{W}_{2345}\,\,(\equiv\mathscr{V})$           &  $\R^4$ bundle over $\CP^1\sqcup \CP^1$ & horizontal foliations                  \\ \hline
  $\mathbf{W}_{12345}$                                 &  $S^1\times S^1$ bundle over $S^3$      & horizontal distributions of type $D_1$ \\ \hline
  $\mathbf{W}_{12456}$                                 &  $S^1\times S^1$ bundle over $S^3$      & vertical distributions of type $D_1$    \\ \hline
  $\mathbf{W}_{123456}$                                &  generic point in  $\Gr_2(\R^6)$        & generic OPS                            \\ \hline
\end{tabular}
\end{teor}\smallbreak

\nit{Proof} We analyse separately the types of OPS's. \smallbreak

\textbf{Type (0,0)}. We first determine the sets
$\mathbf{T}_{00}\cap\mathbf{C}_i:=\mathbf{C}_i^0$.\smallbreak

\noindent --$\mathbf{C}_1^0$. In view of Theorem \ref{opercrit} the
restriction of the map $\tau_h:\mathcal{H}\to \Lambda^2 \mathcal{V}$
must be identically zero, $\Lambda^2 \mathcal{V}$ is generated by
$v_1\wedge v_2$. The choice of basis of $\mathcal{H}$ imposes that
$\nabla v_5$ and $\nabla v_6$ have no skew part, and the projection
of $\nabla v_3$ and $\nabla v_4$ on $\Lambda^2 \mathcal{V}$ produces
a linear combination of $\beta_2$ and $\beta_3$. The lack of this
component means that the projection of $v_1\wedge v_2$ on
$\Lambda^2\mathscr{K}$ should be orthogonal to $\beta_1$ and $\beta
_2$. In conclusion $\mathbf{f}\wedge \mathbf{g}$ must satisfy
equation \rf{g34bt}, so  $\mathbf{C}_1^0$ is determined by the
condition $(\ast)$.\smallbreak

\noindent --$\mathbf{C}_2^0$. The space $S^2_0\mathcal{V}$ is
generated by $v_1\odot v_2$ and $v_1^2-v_2^2$. An OPS is in
$\mathscr{C}_2$ if \begin{equation}(\tau(\mathbf{h}),v_1\odot
v_2)=0=(\tau(\mathbf{h}),v_1^2-v_2^2)\label{van2}\end{equation}\noindent
where $\mathbf{h}$ is a generic element of $\mathcal{H}$. In the
case $(0,0)$ the projection of $\mathcal{H}$ fills the entire space
$\mathscr{K}$, so the equations \rf{van2} must be satisfied by
$\nabla \mathbf{k}$ for any $\mathbf{k}\in \mathscr{K}$.

The computation becomes more explicit if we write the basis
\rf{v1v2} in the
form:$$v_1=\mathbf{f}+ae^5,\,\,\,\,\,v_2=\mathbf{g}+be^5+ce^6.$$\noindent
Keeping in mind Proposition \ref{cov1forms} we are interested in the
"mixed terms" in the space $\mathscr{K}\odot\mathscr{K}^\bot$. The
term
 $v_1^2-v_2^2$ produces $(a\mathbf{f}-b\mathbf{g})\odot e^5-c\mathbf{g}\odot e^6$ and
 $v_1\odot v_2$ gives $(b\mathbf{f}+a\mathbf{g})\odot e^5+c\mathbf{f}\odot e^6$. In these terms
equations \rf{van2} become ($\forall
\mathbf{k}\in\mathscr{K}$):\smallbreak

\noindent $\begin{array}{lcl}
\left\{\begin{array}{l}(J_2\mathbf{k},a\mathbf{f}-b\mathbf{g})-(J_3\mathbf{k},c\mathbf{g})=0
\\
(J_2\mathbf{k},b\mathbf{f}-a\mathbf{g})+(J_3\mathbf{k},c\mathbf{f})=0\end{array}\right.
& \Longrightarrow &
\left\{\begin{array}{l}J_2(a\mathbf{f}-b\mathbf{g})+J_3c\mathbf{g}=0\\J_2(b\mathbf{f}-a\mathbf{g})+J_3c\mathbf{f}=0\end{array}\right.
\end{array}$\smallbreak

\noindent Acting by $J_2$ we
obtain:\begin{equation}\left\{\begin{array}{l}a\mathbf{f}=b\mathbf{g}-cJ_1\mathbf{g}\\a\mathbf{g}=-b\mathbf{f}+cJ_1\mathbf{f}\end{array}\right.\label{eqvan}\end{equation}\noindent
We re-write the first equation as
$a^2\mathbf{f}=b(a\mathbf{g})-c(aJ_1\mathbf{g})$ and substitute the
second and the second transformed by
$J_1$:$$(a^2+b^2-c^2)\mathbf{f}=2bcJ_1\mathbf{f}$$\noindent As
$\mathbf{f}$ is orthogonal to $J_1\mathbf{f}$ a solution is possible
if either $\mathbf{f}=0$ or the coefficients of the left-hand and
the right-hand term both vanish. If $\mathbf{f}=0$, \rf{eqvan}
implies that also $\mathbf{g}=0$ i.e. $v_1\wedge v_2=\pm e^{56}$.
Otherwise must be $a^2+b^2-c^2=0$ and $b=0$ or $c=0$.\smallbreak

If $c=0$ the obvious solution $a=b=c=0$ means that $v_1\wedge v_2\in
\Gr_2(\mathscr{K})$, these are $(2,0)$-type structures which we will
consider later on.\smallbreak

If $b=0$ should be $c=\pm a$, which implies $a\mathbf{g}=\pm
J_1\mathbf{f}$. In other words:$$\mathcal{V}=\langle
\mathbf{f}+ae^5,\,\,\pm(J_1\mathbf{f}+ae^6)\rangle.$$\noindent These
are planes invariant under the action of the standard complex
structure $J_0$ on $\mathscr{I}$ (defined by the $2$-form
$\omega_0=e^{12}+e^{34}+e^{56})$. We know from [8] that
the subset of $\Gr_2(\R^6)$ invariant under the action of $J_0$ is a
disjoint union of two $\mathbb{CP}^2$'s one containing $e^{56}$, the
other $-e^{56}$ and each intersecting $\Gr_2(\mathscr{K})$ in a
$\mathbb{CP}^1$.\smallbreak

\noindent --$\mathbf{C}_3^0$. The space $Tr(\mathcal{V}\otimes
\mathcal{V})$ is generated by $v_1^2+v_2^2$. A computation
completely analogous to the previous case leads to the following
condition determining
$\mathscr{C}_3$:$$a\mathbf{f}=-b\mathbf{g}+cJ_1\mathbf{g}$$
\noindent acting by $J_1$ and summing term by term we get
$\mathbf{g}=-(ab\mathbf{f}+acJ_1\mathbf{f})/(b^2+c^2)$. Denoting by
$d^2=b^2+c^2$ and $\mathbf{f}=d\mathbf{h}$, we see that:$$
\mathcal{V}=\langle
d\mathbf{h}+ae^5,\,\,-a(\frac{b}{d}\mathbf{h}+\frac{c}{d}J_1\mathbf{h})+d(\frac{b}{d}e^5+\frac{c}{d}e^6)\rangle$$\noindent
we can interpret $(b/d,c/d)=(\cos\theta,\sin\theta)$ and thus
re-write this equation
as:\begin{equation}\label{cond3}\mathcal{V}=\langle
d\mathbf{h}+ae^5,\,\,\,e^{J_0\theta}(-a\mathbf{h}+de^5)
\rangle\end{equation} \noindent This expression is easily
interpreted. Given a unit form in $\langle \mathbf{h},e^5\rangle$,
we rotate it in the same plane by $\pi/2$. The resulting form is
rotated by $\theta$ in the direction determined by $J_0$. The choice
of $\mathbf{h}\in \mathscr{K}$, the position in $\langle
\mathbf{h},e^5\rangle$ and the last rotation $\theta$ determine a
point in $\mathbf{C}_3$, which is therefore isomorphic to $S^3\times
S^1\times S^1$.\smallbreak

\noindent --$\mathbf{C}_4^0$. The map
$\tau:V\to\Lambda^2\mathcal{H}$ should be zero. By hypothesis
$\mathbf{f},\mathbf{g},v_5,v_6$ are linearly independent, so their
wedge products generate $\Lambda^2 \mathscr{K}$. The skew part of
$\nabla v_1$ and $\nabla v_2$ is a linear combination of $\beta_2$
and $\beta_3$.  We can avoid skew-terms if $v_1\wedge v_2\in
\Lambda^2\mathscr{K}$. In conclusion $\mathbf{C}_2$ is
$\Gr_2(\mathscr{K})=\mathbf{T}_{20}$ (to be analized
separately).\smallbreak

\noindent --$\mathbf{C}_5^0$. The generators of $S^2_0\mathcal{H}$
with non-zero projection on $\mathscr{K}\odot \mathscr{K}^\bot$ are
$v_3\odot v_4$, $v_3\odot v_5$, $v_3\odot v_6$, $v_4\odot v_5$,
$v_4\odot v_6$. The projections are $a_iv_{iK}\odot e^5+b_i
v_{iK}\odot e^6$, where $a_i$ and $b_i$ denote the $e^5$ and $e^6$
components of $v_3$ and $v_4$. As $J_2v_{1K}\neq J_3v_{1K}$, the
vector $(J_2v_{1K})\odot e^5+ (J_3v_{1K})\odot e^6$ is orthogonal to
the generators (in $\mathscr{K}\odot \mathscr{K}^\bot$) if and only
if $J_2v_{1K}\bot v_{iK}$ and $J_3v_{1K}\bot v_{iK}$ (in
$\mathscr{K}$). As $v_{3K}$, $v_{4K}$, $v_5$, $v_6$ span $\mathscr
K$ we conclude that $\mathbf{C}_5^0=\emptyset$ .\smallbreak

\noindent --$\mathbf{C}_6^0$.  The generator is
$v^2_3+v^2_4+v^2_5+v^2_6$, its  projection on
$\mathscr{K}\odot\mathscr{K}^\bot$ involves $v_3$ and $v_4$. The
projection of both $\nabla v_1$ and $\nabla v_2$ is zero if and only
if $\langle v_5,v_6 \rangle=\langle J_2\mathbf{f},J_3\mathbf{f}
\rangle=\langle J_2\mathbf{g},J_3\mathbf{g}\rangle$.This means that
the OPS's satisfies condition $(\ast)$.\smallbreak

Some relevant intersections and inclusions of the sets
$\mathbf{C}^0_i$ are listed
below:$$\mathbf{C}^0_1\cap\mathbf{C}^0_2=\CP^2\sqcup\CP^2,\,\mathbf{C}^0_1\cap\mathbf{C}^0_4=\mathscr{D},\,\mathbf{C}^0_3\subset\mathbf{C}^0_1\subset\mathbf{C}^0_6,$$
$$\mathbf{C}^0_4\subset\mathbf{C}^0_2\subset\mathbf{C}^0_6,\mathbf{C}^0_2\cap\mathbf{C}^0_3=\mathbf{C}^0_3\cap\mathbf{C}^0_4=\emptyset.$$\smallbreak\noindent The higher-order intersections of
$\mathbf{C}_1^0$ can be easily deduced. The remaining types of
structures are analized directly in terms of $\mathscr{W}_i$ instead
of $\mathscr{C}_i$. The reader should keep in mind the bases given
in Table 1, Theorem \ref{opercrit} and Proposition
\ref{cov1forms}.\smallbreak

\noindent \textbf{Type (0,2)}. We see immediately that these OPS's
belong (and exhaust) $\mathbf{W}_{4}$. \smallbreak

\noindent \textbf{Type (2,0)}.  The lack of ``mixed'' terms
determines the absence of $\mathscr{W}_3$ and $\mathscr{W}_6$
components. As $\mathscr{K}^\bot\cap\mathcal{V}=\emptyset$ there are
no $\mathscr{W}_2$ and $\mathscr{W}_4$ components. There is always a
$\mathscr{W}_5$ component and $\mathscr{W}_1$ vanishes for
$J_1$-invariant planes (in $\mathscr{D}$). \smallbreak

\noindent \textbf{Type (1,1)}. The absence of $\mathscr{W}_1$ is
obvious. The lack of ``mixed'' terms determines again the absence of
$\mathscr{W}_3$ and $\mathscr{W}_6$. Observe
that:\begin{equation}\label{jjjj}
\begin{array}{l}\nabla (J_1\mathbf{f})=-(J_3\mathbf{f})\odot e^5+(J_2\mathbf{f})\odot e^6\\ \nabla
(J_2\mathbf{f})=-\mathbf{f}\odot e^5-(J_1\mathbf{f})\odot e^6\\
\nabla (J_3\mathbf{f})=(J_1\mathbf{f})\odot e^5-\mathbf{f}\odot
e^6.\end{array}\end{equation} \noindent Comparing the above terms to
$\Lambda^2\langle J_1\mathbf{f},J_2\mathbf{f},J_3\mathbf{f}\rangle$,
an easy computation shows that $\mathscr{W}_2$ and $\mathscr{W}_5$
are always present. In conclusion the $(1,1)$ structures belong to
$\mathscr{W}_{245}$.\smallbreak

\noindent \textbf{Type (0,1)}. The lack of $\mathscr{W}_1$ component
is obvious. Comparing $\mathbf{f}\odot \mathbf{d}$ and
$\mathbf{f}\odot \mathbf{e}$ to \rf{jjjj} we find always a
$\mathscr{W}_2$ and $\mathscr{W}_3$ components. The presence of
$\mathscr{W}_5$ is obvious and the presence of $\mathscr{W}_4$ is
motivated in a way analogous to the one used in the previous case.
The $\mathscr{W}_6$ component is absent as the relevant term in
$Tr(\mathcal{H}\otimes \mathcal{H})$ is $\mathbf{f}\odot \mathbf{e}$
which is always orthogonal to $\nabla\mathbf{f}$. In conclusion the
$(0,1)$ structures belong to $\mathscr{W}_{2345}$.\smallbreak

\noindent \textbf{Type (1,0)}. The set
$\mathbf{T}_{10}\cap\mathscr{C}_1$ is defined by condition $(\ast)$.
The relevant terms in $S^2_0\mathcal{V}$ are $\mathbf{f}\odot
\mathbf{e}$ and $-\mathbf{g}\odot \mathbf{e}$. The projection of
$\mathcal{H}$ onto $\mathscr{K}$ is $\langle
\mathbf{f}\rangle^\bot=\langle
J_1\mathbf{f},J_2\mathbf{f},J_3\mathbf{f}\rangle$, so the $(1,0)$
structures do have $\mathscr{W}_2$ component. Analogous
considerations show that there are $\mathscr{W}_4$ and
$\mathscr{W}_5$ components.

The term generating $\mathscr{W}_6$ is $\mathbf{g}\odot \mathbf{e}$,
orthogonal to  $\nabla \mathbf{g}$. Assuming condition $(\ast)$,
Proposition \ref{cov1forms} implies that there is no trace. If
$\mathbf{g}\in \langle J_2\mathbf{f}, J_3\mathbf{f}\rangle$ we can
write \begin{equation}\mathcal V=\langle \mathbf{f}\, \,,\,\,
Y_1J_2\mathbf{f}+Y_2J_3\mathbf{f}+Y_3e^5+Y_4e^6\rangle \,
,\,\,\,\,\,\,\,\,\,\sum Y_i^2=1\label{gjj}\end{equation}\noindent
which, up to a proportionality factor, means
\begin{equation}f=Y_1J_2\mathbf{g}+Y_2J_3\mathbf{g}\label{fjjj},\end{equation}$$\nabla
\mathbf{f}=-Y_1\mathbf{g}\odot e^5-Y_2\mathbf{g}\odot
e^6+Y_1J_1\mathbf{g}\odot e^5-Y_2J_1\mathbf{g}\odot e^6.$$ \noindent
Compare $-(Y_1\mathbf{g}\odot e^5+Y_2\mathbf{g}\odot e^6)$ to $Y_3
\mathbf{g}\odot e^5+Y_4 \mathbf{g}\odot e^6$. The trace vanishes
if\begin{equation}\label{perp}Y_1 Y_3+Y_2Y_4=0.\end{equation}The
part in $\mathscr{W}_3$ is generated by $\mathbf{g}\odot
\mathbf{e}$, so $\nabla \mathbf{g}$ gives no contribution. The
non-zero contribution comes from $\nabla v_3$ and $\nabla v_4$. If
$\mathbf{f}=J_1\mathbf{g}$, we can set $v_3=J_2\mathbf{g}$,
$v_4=J_3\mathbf{g}$ and \rf{jjjj} implies that there is trace.
Setting $\mathbf{f}=Y_1J_2\mathbf{g}+Y_2J_3\mathbf{g}$,
$v_3=J_1\mathbf{g}$ and $v_4=Y_2J_2\mathbf{g}-Y_1J_3\mathbf{g}$, we
see that $\nabla v_3$ gives no contribution and we just compare
$Y_3\mathbf{g}\odot e^5+Y_4 \mathbf{g}\odot e^6$ to
$Y_2\mathbf{g}\odot e^5-Y_1\mathbf{g}\odot e^6$ in $-\nabla v_4$.
The trace vanishes if \begin{equation}\label{par}Y_2 Y_3-Y_1Y_4=0
.\end{equation} \noindent So if we take $\mathbf{f}\in S^3$, a
$(1,0)$ structure with $\mathbf{g}$ generically in $\langle
J_1\mathbf{f},J_2\mathbf{f},J_3\mathbf{f}\rangle$ belongs to
$\mathscr{W}_{123456}$. Condition $(\ast)$ detects a $\CP^1\sqcup
\CP^1$ of planes in $\mathscr{K}$ and no condition on $\mathbf{e}$,
so $\mathbf{T}_{10}\cap \mathbf{W}_{2345}$ is an $\R^2$ bundle over
$\CP^1\sqcup \CP^1$ (a subset of $\mathscr{V}$).\smallbreak

\noindent Condition \rf{gjj} determines in $\mathbf{T}_{10}$ a
subbundle over $S^3$ with fibre $S^3$. Equations \rf{perp} and
\rf{par} define cones in the same space. So the fibres are
intersections of quadrics in $\R^4$, a sphere and a cone. An easy
computation shows that with a suitable choice of basis the
intersection of \rf{perp} with the sphere is defined
by:$$X_1^2+X_2^2=\frac{1}{2},\,\,\,\,\,X_3^2+X_4^2=\frac{1}{2}.$$
\noindent We conclude that the sets $\mathbf{T}_{10}\cap
\mathbf{W}_{12345}$ and $\mathbf{T}_{10}\cap \mathbf{W}_{12456}$ are
both diffeomorphic to an $S^1\times S^1$ bundle over $S^3$. Observe
that $Y_1^2+Y^2_2=0$ and $Y_3^2+Y^2_4=0$ belong to the intersection
of \rf{gjj} and \rf{perp} but also \rf{gjj} and \rf{par}. These are
$(1,1)$ and $(2,0)$-type structures respectively. So on the
intersection of the sets $\mathbf{W}_{12345}$ and
$\mathbf{W}_{12456}$ both classes ``decay'' to $\mathbf{W}_{245}$
and $\mathbf{W}_{15}$.\smallbreak

We complete the proof highlighting the following facts. The elements
of $\mathbf{C}_3$ (recall \rf{cond3}) satisfy condition $(\ast)$, so
$\mathbf{C}_3=\mathbf{C}_{136}$. Setting $a=0$ in \rf{cond3} we
obtain a generic element in $\mathbf{T}_{11}$, so this condition
detects completely $\mathbf{W}_{245}\cong S^3\times S^1\times S^1$
including $\mathbf{T}_{11}\cong S^3\times S^1$. The $J_0$ invariance
is a particular case of condition $(\ast)$ so
$\mathbf{C}^0_{126}=\mathbf{W}_{345}$. In view of Remark \ref{type}
$\mathbf{W}_{2345}$ coincides with $\mathscr{V}$. \smallbreak

The column ``geometric properties'' is given just for completeness,
for more details see [11]. Observe that $\mathbf{W}_{2345}$
and its subvarieties contain horizontal foliations and
$\mathbf{W}_{15}$ contains vertical foliations. \qed\smallbreak

\noindent The compactification leads a null-torsion class to
``decay'' to a ``smaller''-one. In particular $\mathscr{W}_{2345}$
falls into $\mathscr{W}_{345}$, $\mathscr{W}_{245}$ ect. on disjoint
intrinsic torsion sub-varieties. The same occurs with
$\mathscr{W}_{5}\subset\mathscr{W}_{15}$. Similarly
$\mathbf{T}_{10}$ is compactified by $\mathbf{T}_{20}$ and
$\mathbf{T}_{11}$.

\section{Constant torsion varieties}

We now aim to illustrate the property of ITV's claimed in Remark
\ref{ctvrem}. Given an OPS $\mathfrak{P}_1$ any other OPS
$\mathfrak{P}_2$ can be obtained by a base change in $O(N)$ such
that $\mathcal V_1\longrightarrow \mathcal{V}_2$ and
$\mathcal{H}_1\longrightarrow \mathcal{H}_2$. Keeping this action in
mind we combine Equation \rf{ctveq} with Theorem \ref{opercrit} to
get the following useful criterion:

\begin{corol}\label{opcrit2}Given an OPS $\mathfrak P$ and $g\in O(N)$, $g\mathfrak P$ belongs to the same $CTV$ if and only if for each vertical or horizontal 1-form $v$ occurs that $\nabla g v=g(\nabla v)$.\end{corol}

\noindent So in our case we are looking for transformations $g\in
SO(6)$ satisfying the above conditions. The tangent space
$T\mathscr{I}$ has an invariant basis $\{e_i\}$ compatible with the
natural splitting $\mathscr{K}\oplus \mathscr{K}^\bot$ (Malchev
basis) \emph{and} with the standard complex structure on
$\mathscr{I}$ given by the complex Heisenberg group. The Levi-Civita
connection involves a tern ($J_1$,$J_2$,$J_3$) of orthogonal almost
complex structures on $\mathscr{K}$ which satisfy the quaternion
product rule. In other words, there is a standard $Sp(1)$ structure
on $\mathscr{K}$. Using the notation of \rf{v1v2}, a simple
straightforward computation leads to:

\begin{lemma} Given a 1-form $v$ on $\mathscr{I}$ and $\mathbf{k}\in\mathscr K$, we define $g\in SO(6)$ by:$$\begin{array}{lcllcl}
g\mathbf{f}&=&\cos\theta \mathbf{g}+\sin\theta J_1\mathbf{k},\,\,\,\,\,\, & gJ_1\mathbf{f}&=& -\sin\theta \mathbf{k}+\cos\theta J_1\mathbf{k}, \\ gJ_2\mathbf{f} &=& \cos\theta J_2\mathbf{k}+\sin\theta J_3\mathbf{k},\,\,\,\,\,\, & gJ_3\mathbf{f}&=& -\sin\theta J_2\mathbf{k}+\cos\theta J_3\mathbf{k}, \\
ge^5&=&\cos\theta e^5+\sin\theta e^6,\,\,\,\,\,\,  & ge^6&=&
-\sin\theta e^5+\cos\theta e^6 \end{array}.$$ Then $\nabla v$ keeps
its form under this ``coupled'' action of $Sp(1)\times SO(2)$
i.e.\smallbreak

$g((J_2\mathbf{f})\odot e^5+(J_3\mathbf{f})\odot
e^6+a\beta_2+b\beta_3)=(J_2\mathbf{k})\odot e^5+(J_3\mathbf{k})\odot
e^6+a\beta_2+b\beta_3.$\end{lemma}

The requirement of Corollary \ref{opcrit2} is satisfied if and only
if $\theta=0$ so:

\begin{prop} The CTV's of $SO(4)\times SO(2)$-structures on $\mathscr{I}$ are $Sp(1)$-orbits. \end{prop}

\noindent Considering the specific compatible bases used in the
previous sections, the components in $\mathscr{K}$ undergo the
$Sp(1)$-action and the remaining parameters determine the
cross-section in the set of $SO(6)$-orbits (Remark \ref{ctvrem}). Up
to $\mathbb{Z}^2$ overabundance:\smallbreak

\noindent $\mathbf{W}_4$ is a degenerate CTV as $Sp(1)\times Id$ is
contained in the reduced structure group.

\noindent $\mathbf{W}_{5}$ is composed by two $\CP^1$-s, each a
constant torsion subvariety of $\mathbf{W}_{15}$.

\noindent $\mathbf{W}_{15}$ is a $\CP^1$ bundle (the coadjoint orbit
of $Sp(1)\cong SU(2)$) over $\CP^1$.

\noindent $\mathbf{W}_{345}$ is an $S^3$ bundle (the orbit of
$\mathbf{f}$) over $S^1$ (the choice of $\mathbf{d}$ in \rf{v1v2}).

\noindent $\mathbf{W}_{245}$ is a trivial $S^3$  bundle (the orbit
of $\mathbf h$ in \rf{cond3}) over $S^1\times S^1$.

\noindent $\mathbf{W}_{2345}$ is a $\CP^1$ bundle over $\R^4$
compactified as in Remark \ref{type}.

\noindent $\mathbf{W}_{12456}$ - the $Sp(1)$-orbit of $\mathbf{f}$
is an $S^3$ the remaining parameters determine $S^1\times S^1$.

\noindent $\mathbf{W}_{12345}$ is analogous to the previous case.

\noindent $\mathbf{W}_{123456}$ is a $\CP^1$-bundle over $\CP^1$
(the set of orbits in $\mathscr{K}$) times compactified $\R^4$.

\section{Toric invariance and moment mapping}

Given a maximal torus $T\subset O(2N)$, a $T$-structure is a
splitting of the tangent space into orthogonal real two-dimensional
subspaces. The standard $Sp(1)\times SO(2)$-structure on
$\mathscr{I}$ determines a natural invariant $T^3$-structure (say
$\bar{T}$): $$P_1=\langle e^1,e^2\rangle,\,\,\,\, P_2=\langle
e^3,e^4\rangle,\,\,\,\, P_3=\langle e^5,e^6\rangle,$$ \noindent The
Levi-Civita derivative of Proposition \ref{cov1forms} is resumed by
maps:$$\nabla:P_1\longrightarrow P_2\odot P_3,
\,\,\,\,\,\,\,\,\nabla:P_2\longrightarrow P_1\odot P_3,
\,\,\,\,\,\,\,\nabla:P_3\longrightarrow P_1\wedge P_2$$\noindent
This fact motivates the following:\begin{teor}\label{torinvar}The
ITV's of invariant $SO(4)\times SO(2)$-structures on $\mathscr{I}$
are stable under the action  of $\bar{T}$ with the exception of
$\mathbf{W}_{12345}$ and $\mathbf{W}_{12456}$.
\end{teor}\nit{Proof} The action of $\bar {T}$ preserves the splitting
$\mathscr{K}\oplus\mathscr{K}^\bot$ and the type of structure. Some
conditions defining ITV's involve $J_0$ (recall \rf{cond3} ) and
$J_i$-s, but they commute as elements as of $\bar T$. This proves
the invariance of $\mathbf{W}_{245}$ and  $\mathbf{W}_{4}$,
$\mathbf{W}_{5}$, $\mathbf{W}_{15}$, $\mathbf{W}_{345}$,
$\mathbf{W}_{2345}$. Analogously we see that $\mathbf{W}_{12345}$
and $\mathbf{W}_{12456}$ are stable under the action of the subtorus
of $\bar T$ acting on $\mathscr{K}$. On the contrary the action of a
generic element of $\bar T$ on the space $\mathscr{K}^\bot$ breaks
the ``coupling'' \rf{perp} and \rf{par}.\qed\smallbreak

As a consequence the ``geography'' of $\Gr_2(\R^6)$ can be neatly
visualized.\smallbreak

The classes of Hermitian structures on six-dimensional nilmanifolds
has been described in [1] in terms of faces, edges and
other segments in a solid tetrahedron. An accurate geometric
interpretation of this fact has been given in [8, 10].
The reductions of a Riemannian structure determined by a $2$-form
are parametrized by $SO(N)$ (co)adjoint orbits. The $2$-form in
question can be seen as the image of a standard moment mapping
associated to the Konstant--Kirillov--Souriau symplectic structure.
Typical constructions in symplectic geometry, as symplectic
fibrations, quotients etc. can be interpreted in terms of
compatibility of $G$-structures, and appear as efficient tools for
their parametrization. The restriction of the structural group
$SO(2N)$ to a maximum torus gives rise to a Hamiltonian toric action
with moment map $\mu_T$ to the dual of the Lie algebra
$\mathfrak{t}^*$. The celebrated Atiyah and Guillemmin--Sternberg
Convexity Theorem, implies that the image of each classifying orbit
by $\mu_T$ is a specific polytope.\smallbreak

\begin{figure}[!h]
\centering\includegraphics[width=0.7\textwidth]{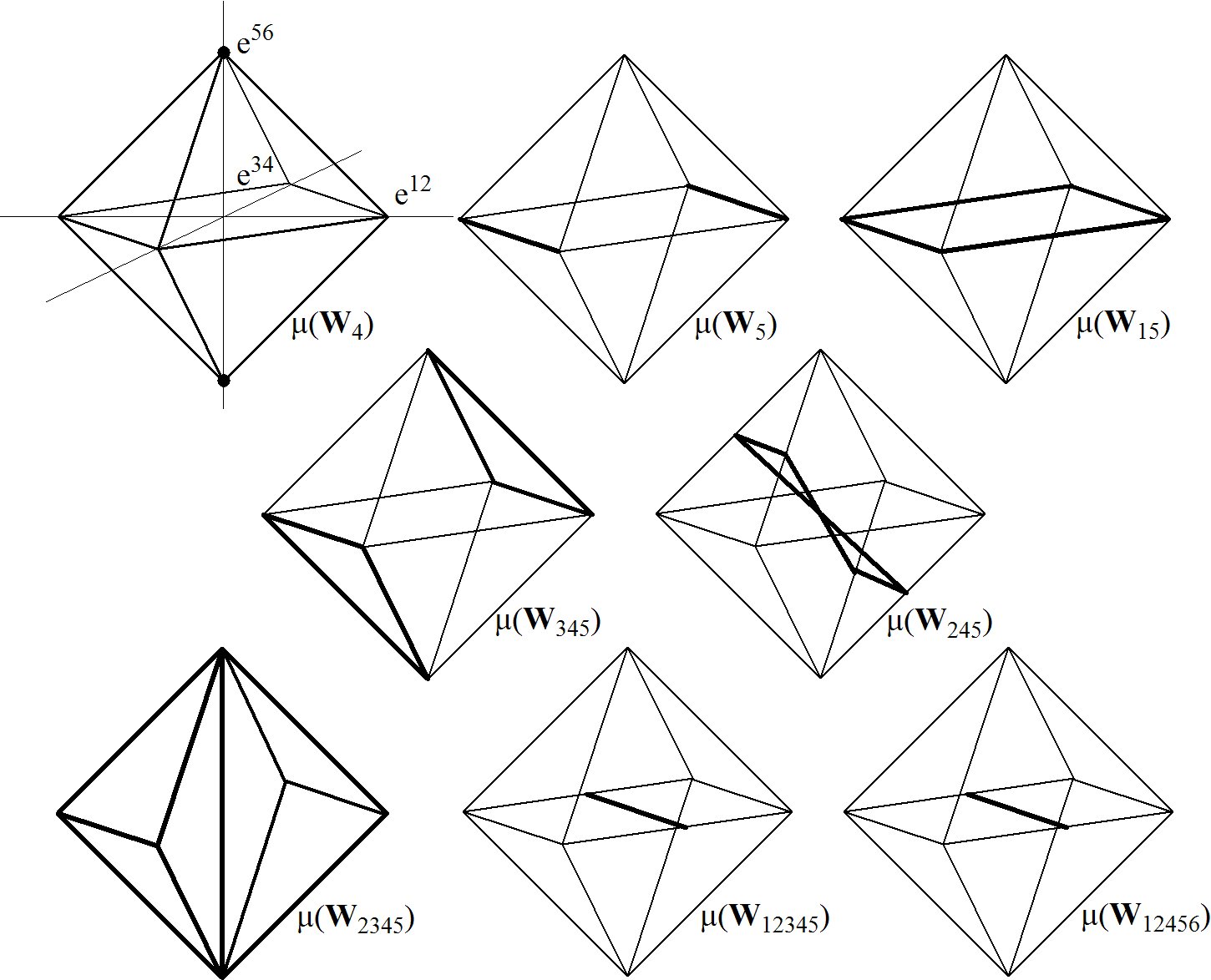}\\
\caption{Images by $\mu_T$ of the ITV's of OPS'S in the moment
polytope of $\Gr_2(\R^6)$.} \label{torsgras}\end{figure}

In our case the moment mapping associated to $\bar{T}$ is the
orthogonal projection onto the subspace
$\Lambda^2(\R^6)\supset\langle
e^{12},e^{34},e^{56}\rangle\cong\mathfrak{t}^*\cong \R^3$. The image
of $\Gr_2(\R^6)$ is an octahedron with vertices $\pm e^{12},\pm
e^{34},\pm e^{56}$ (see Fig. \ref{torsgras}). As the dimension of
$\Gr_2(\R^6)$ is higher than twice the dimension of $\bar{T}$, the
Hamiltonian $T$-action gives rise to a non-trivial symplectic
quotient $\mathcal{Q}_P=\mu_T^{-1}(P)/\bar T $ at each point $P$ of
the octahedron. In our case a generic $\mathcal{Q}_P$ is
diffeomorphic to a real $2$-sphere (see [10]). By analogy,
in [9] we have defined the \emph{intrinsic quotient}
$\mathcal{IQ}_P$ of an ITV to be the subset of $\mathcal{Q}_P$ which
represents the ITV. The construction (parametrizing the set of
$T$-orbits in the inverse image of a point in the moment polytope)
is valid also if the ITV is not a symplectic manifold.

\begin{teor} For the ITV's of invariant $SO(2)\times SO(4)$-structures on $\mathscr{I}$ occurs:\smallbreak
\begin{tabular}{|l|l|l|}
  \hline ITV                    &  $\mu_T(\mathbf{W}_{ijk..})$ (see Fig. \ref{torsgras})                           & $\mathcal{IQ}_P$ \\ \hline \hline
  $\mathbf{W}_4$                &  two points                                                                      & a point        \\ \hline
  $\mathbf{W}_5$                &  two disjoint line segments                                                      & a point        \\ \hline
  $\mathbf{W}_{15}$             &  filled square                                                                   & a point        \\ \hline
  $\mathbf{W}_{345}$            &  two disjoint triangular faces                                                   & a point        \\ \hline
  $\mathbf{W}_{245}$            &  two triangles in the plane $x+y+z=0$                                            & a point        \\ \hline
  $\mathbf{W}_{2345}$           &  two tetrahedra with common edge                                                 & a point        \\ \hline
  $\mathbf{W}_{12345}$          &  line segment                                                                    & a curve \\ \hline
  $\mathbf{W}_{12456}$          &  line segment                                                                    & a curve \\ \hline
  $\mathbf{W}_{123456}$         &  the octahedron                                                                  & $S^2$          \\ \hline
\end{tabular}
\end{teor}

\smallbreak

\nit{Proof}  The images of each ITV are obtained by elementary
computations taking the projection on $\langle
e^{12},e^{34},e^{56}\rangle$ of the simple $2$-form $v^1\wedge v^2$
associated to the OPS. In particular the points of $\mathbf{W}_{4}$
are fixed by the action of $\bar{T}$, $\mathbf{W}_{5}$,
$\mathbf{W}_{15}$ and $\mathbf{W}_{345}$ are toric varieties, and as
such can be represented as (well known) $\bar{T}$-orbit fibrations
over their images by $\mu_T$. As $\mathbf{W}_{245}$ and
$\mathbf{W}_{2345}$ are $\bar{T}$-stable and their intrinsic
quotient is determined just by a dimensional check.

The sets $\mathbf{W}_{12345}$ $\mathbf{W}_{12456}$ are not $\bar
T$-stable, but as $J_2$ and $J_3$ commute with the restriction $\bar
T^2$ of $\bar T$ to $\mathscr{K}$, equation \rf{fjjj} is $\bar
T^2$-invariant. The $\bar T^2$-orbit of the unit form $f=\sum_1^4
x_ie^i$ depends on the value of $x_1^2+x_2^2$ in $[0,1]$. The
generic $\bar T^2$-orbit in the base space of the fibration ($S^3$)
is a two-torus. Values 0 and 1 correspond to degenerate $S^1$ orbits
in $\langle e^3, e^4\rangle$ and $\langle e^1, e^2\rangle$. The
$\bar T^2$-orbits of $Y_1J_2\mathbf{f}+Y_2J_3\mathbf{f}\in
\mathscr{K}$ are analogously described. The couple $(Y_1,Y_2)$
determines distinguished $\bar T^2$-orbits on the fibres but it is
easy to see, via standard parametrization of the $\bar T^2$-action
on $\mathscr{K}$, that the orbit depends on $Y_1^2+Y_2^2\in (0,1)$.
The set of orbits is parametrized by $[0,1]\times(0,1)$. Both
$\mathbf{W}_{12345}$ and $\mathbf{W}_{12456}$ project on $\langle
e^{12},e^{34}\rangle$ onto the line:\begin{equation}\label{proj}
X=Y_1(x_1x_4+x_2x_3)+Y_2(x_2x_4-x_1x_3)=-Y.\end{equation}

\noindent   The projection imposes a one-dimensional constrain, so
the set of orbits in the inverse image of each point of
$\mu_T(\mathbf{W}_{12345})$ and $\mu_T(\mathbf{W}_{12456})$ is
one-dimensional.\qed\smallbreak

For a more precise description of the intrinsic quotient of
$\mathbf{W}_{12345}$ and $\mathbf{W}_{12456}$, an explicit
parametrization of the symplectic quotient of $\Gr_2(\R^6)$ is
needed.

\begin{rem} The inclusions of ITV's are very well represented by the moment polytopes. The compactification described in Remark \ref{type} can be also  visualized, $\mu_T(\mathbf{T}_{11})$ is the origin, $\mu_T(\mathbf{T}_{01})$ is the vertical segment connecting $\pm e^{56}$ etc.
\end{rem}

The author believes that the rich geometry, which emerges in the
analysed case, will stimulate the construction of a general theory
on Intrinsic Torsion Varieties.

\section*{Acknowledgement}
The author thanks Simon Salamon for his help.

\section*{References}
\small
\noindent $[1]$ E. Abbena, S. Garbiero, and S. Salamon. Almost Hermitian geometry of 6-dimensional
nilmanifolds. \emph{Ann. Sc. Norm.} Sup., $30,\,147 - 170, 2001$.\medbreak

\noindent $[2]$ D. E. Blair. Geometry of manifolds with structural group $U(n)\times O(s)$. \emph{Journal of
differential geometry}, $4,\, 155 - 167,\, 1970$.\bigbreak

\noindent $[3]$ R. Bryant. Calibrated embeddings in the special lagrangian and coassociative cases.
\emph{Ann. Global Anal. Geom. -- Special issue in memory of Alfred Gray (1939-1998)}, $18$, no.
$3-4,\,405-435,\, 2000$.\bigbreak

\noindent $[4]$ F. M. Cabrera and A. Swann. The intrinsic trosion of almost quaternion-hermitian
manifolds. \emph{Syddansk Universitet}, Preprint No.6 ISSN No. $0903-3920$, July $2007$.\bigbreak

\noindent $[5]$ S. Chiossi and A. Swann. $G_2$ -structures with torsion from half-integrable nilmanifolds.
\emph{Journal of Geometry and Physics}, $54,\,265-285,\, 2005$.\bigbreak

\noindent $[6]$ D. Conti-A.Tomassini. Special symplectic six-manifolds. \emph{Quart J. Math}, $58,\,297-311,\,2007$.\bigbreak

\noindent $[7]$ A. Gray and L. Hervella. The sixteen classes of almost Hermitian manifolds and their
linear invariants. \emph{Ann. Mat. Pura Appl.}, $123,\,35-58,\, 1980$.\bigbreak

\noindent $[8]$ G. Mihaylov. Special Riemannian structures in six dimensions. \emph{Doctoral thesis, University
of Milan}, $2008$.\bigbreak

\noindent $[9]$ G. Mihaylov. Toric moment mappings and Riemannian structures. \emph{Geometriae Dedicata},
ISSN $0046-5755$, DOI $10.1007/s10711-012-9720-6, 2012$.\bigbreak

\noindent $[10]$ G. Mihaylov and S. Salamon. Intrinsic torsion varieries. \emph{Note di matematica Universit\'a di Lecce}, Suppl. N. $1,\,349-376.,\, 2009$.\bigbreak

\noindent $[11]$ A. M. Naveira. A classification of Riemannian almost-product manifolds. \emph{Rendiconti
di Matematica Roma}, $3,\,577-592, 1983$.\bigbreak

\noindent $[12]$ S. Salamon. A tour of exceptional geometry. \emph{Milan J. Math.}, $59-94,\,311-333,\, 2003$.\bigbreak

\noindent $[13]$ S. Sternberg. Lectures on Differential Geometry. \emph{New York Chelsea Publishing Company},
$1983$.\bigbreak

\noindent $[14]$ K. Yano. On a structure defined by a tensor field f of type $(1,1)$ satisfying $f^3+f = 0$. \emph{Tensor}, $14,\,99-109,\, 1963$.

\end{document}